\documentclass[12pt,a4paper]{amsart}
\usepackage{amssymb,amscd,amsxtra,calc}

\theoremstyle{plain}
    \newtheorem{thm}{Theorem}[section]
    \newtheorem{claim}[thm]{Claim}
    \newtheorem{corollary}[thm]{Corollary}

    \newtheorem{lemma}[thm]{Lemma}

    \newtheorem{theorem}[thm]{Theorem}

\theoremstyle{definition}

\theoremstyle{remark}

\newcommand{\bA}{\mathbb{A}}

\newcommand{\bG}{\mathbb{G}}

\newcommand{\PP}{\mathbb{P}}
\newcommand{\Q}{\mathbb{Q}}

\newcommand{\Aut}{\operatorname{Aut}}

\newcommand{\Lie}{\operatorname{Lie}}

\newcommand{\op}{\operatorname{op}}

\begin{document}

\title[Log Kodaira dimension of homogeneous varieties]
{Log Kodaira dimension of homogeneous varieties}

\author{Michel Brion}
\address
{
\textsc{Institut Fourier, Universit\'e de Grenoble,}\endgraf
\textsc{
B.P. 74, 38402 Saint-Martin d'H\`eres, France
}}
\email{Michel.Brion@ujf-grenoble.fr}

\author{De-Qi Zhang}
\address
{
\textsc{Department of Mathematics, National University of Singapore,} \endgraf
\textsc{
10 Lower Kent Ridge Road,
Singapore 119076
}}
\email{matzdq@nus.edu.sg}

\begin{abstract}
Let $V$ be a complex algebraic variety, homogeneous under the action
of a complex algebraic group. We show that the log Kodaira dimension
of $V$ is non-negative if and only if $V$ is a semi-abelian variety.
\end{abstract}


\thanks{The second-named author is supported by an ARF of NUS}

\maketitle

\section{Introduction}\label{Intro}
Throughout this short note, we work over the field of complex numbers.
Given a smooth quasi-projective variety $V$,
we recall the definition of the log Kodaira dimension
$\overline{\kappa}(V)$. Take any projective variety $X$ which contains
$V$ as a dense open subset. By blowing up subvarieties of the (reduced) boundary
$D := X \setminus V$, we may assume that $X$ is smooth and $D$ is
a simple normal crossing divisor. In this case, the pair $(X, D)$
is called a {\it smooth projective compactification of} $V$.
Let $K_X$ be a canonical divisor of $X$. The {\it log Kodaira dimension}
$\overline{\kappa} = \overline{\kappa}(X, K_X + D)$ is defined
as the unique value
\[ \overline{\kappa}  \in \{ -\infty, 0, 1, \dots, \dim V \} \]
such that for some positive constants $\alpha, \beta$, we have
\[
\alpha m^{\overline{\kappa}} \le h^0(X, m(K_X + D)) \le \beta m^{\overline{\kappa}}
\]
for any sufficiently large and divisible positive integer $m$.
As $\overline{\kappa}$ is independent of the choice of the smooth
projective compactification $(X, D)$ (cf.~\cite[\S 11.1]{Ii81}), we may
set $\overline{\kappa} := \overline{\kappa}(V)$.

Also, recall that a Cartier divisor $D$ on a projective variety $X$
is {\it pseudo-effective} if its numerical equivalence class is contained
in the closure of the convex cone spanned by the effective divisor
classes on $X$.

We may now state the main result of this note:

\begin{theorem}\label{ThA}
Let $X$ be a smooth projective variety, and $G$ a connected algebraic
group (possibly nonlinear) of automorphisms of $X$. Assume that
$G$ has a dense open orbit $V$ in $X$, and $D := X \setminus V$
is a simple normal crossing divisor.
Then the following conditions are equivalent.
\begin{itemize}
\item[(1)] The log canonical divisor $K_X + D$ is pseudo-effective.
\item[(2)] $K_X + D$ is linearly equivalent to $0$.
\item[(3)] The log Kodaira dimension $\overline{\kappa}(V)$ is non-negative.
\item[(4)] $G$ is a semi-abelian variety.
\end{itemize}
\end{theorem}

We recall that a {\it semi-abelian variety} is an algebraic group $G$
that lies in an exact sequence
\[ 1 \longrightarrow T \longrightarrow G \longrightarrow A \longrightarrow 1, \]
where $T$ is an algebraic torus (isomorphic to a product of copies
of the multiplicative group $\bG_m$), and $A$ an abelian variety;
then $G$ is connected and commutative.

Here are a few words about the proof of Theorem \ref{ThA}, which is very simple.
If $V$ is not a semi-abelian variety, we show that the general member of some
covering family of affine lines on $G$ is not contractible to a point in $V$.
Chevalley's structure theorem for algebraic groups is also used, so our proof
is quite geometric.

From Theorem \ref{ThA}, we derive the following purely group-theoretic
statement:

\begin{corollary}\label{CorA}
Let $G$ be a connected algebraic group, and $H \le G$ a closed subgroup.
Then $\overline{\kappa}(G/H) \le 0$, with equality if and only if $H$ contains
a closed normal subgroup $N$ of $G$ such that $G/N$ is a semi-abelian variety.
\end{corollary}

\section{Proofs of Theorem \ref{ThA} and Corollary \ref{CorA}}

The following result is well-known. See \cite[Lemma 5.11]{KeMc99}
for a generalization to singular pairs.

\begin{lemma}\label{KeMc}
Let $Y$ be a smooth projective variety and $A$ a simple normal crossing
divisor on $Y$. Suppose that $C$ is the general member of a covering family
of rational curves on $Y$. Let $d$ be the number of times that $C$ meets
$A$, i.e., the degree of the (reduced) support of the divisor $\nu^*(A_{|C})$, where
$\nu : \PP^1 \to C$ is the normalization. Then we have:
\begin{itemize}
\item[(1)]
If $d \le 2$, then $(K_Y + A) \cdot C \le 0$.
\item[(2)]
If $d \le 1$, then $(K_Y + A) \cdot C < 0$. In particular $K_Y + A$ is not
pseudo-effective.
\end{itemize}
\end{lemma}


{\it We now prove Theorem \ref{ThA}.}

(2) $\Rightarrow$ (3) and
(3) $\Rightarrow$ (1) follow readily from the definitions of the log
Kodaira dimension and pseudo-effectivity.

(4) $\Rightarrow$ (2). Since the semi-abelian variety $G \subseteq \Aut(X)$ is commutative
and acts faithfully on $X$,
the stabilizer of any point of $V$ is trivial, and hence $V \cong G$.
Also, the action of the Lie algebra of $G$ on $X$ yields a trivialization
of the log tangent bundle $T_X( - \log D)$ (see \cite[Main Thm.]{Wi04}
and \cite[Thm.~2.5.1]{Br07}). In particular, the determinant of
$T_X(- \log D)$ is trivial as well, i.e., $K_X + D \sim 0$.

(1) $\Rightarrow$ (4).
Fix a point $x_0 \in V$ and denote by $H \le G$ its isotropy subgroup,
so that $V$ is $G$-equivariantly isomorphic to the homogeneous space
$G/H$, with both the $G$-actions from the left. Since $G$
acts faithfully on $V$, we have:

\begin{claim}\label{c1}
$H$ contains no non-trivial closed normal subgroup of $G$.
\end{claim}

Next we show:

\begin{claim}\label{c2}
$H$ contains the image of every additive one-parameter subgroup
of $G$.
\end{claim}

We prove Claim \ref{c2}.
Assume the contrary that there exists an additive one-parameter subgroup of $G$,
i.e., a homomorphism of algebraic groups $u : \bG_a \to G$
(where $\bG_a$ denotes the additive group), with image not contained
in $H$. Then the $\bG_a$-orbit $u(\bG_a) \cdot x_0$ is a curve and isomorphic
to $\bG_a/F$, where the isotropy subgroup $F$ is a finite subgroup of
$\bG_a$ and hence is trivial. Thus we obtain an embedding
\[
f : \bA^1 \cong \bG_a \cong u(\bG_a) \cdot x_0 \hookrightarrow G/H, \]
\[ t \longmapsto u(t) \cdot x \]
and hence a rational curve $C \subset X$
which intersects the boundary $D$ at at most one point. Also, the translates of
$C$ by $G$ form a covering family of rational curves on $X$. By Lemma
\ref{KeMc}, it follows that $K_X + D$ is not pseudo-effective,
contradicting our assumption. This proves Claim \ref{c2}.

\par \vskip 1pc
We return to the proof of Theorem \ref{ThA}.
By Chevalley's structure theorem
(see \cite[Thm.~16]{Ro56}), $G$ lies in a unique extension
\[ 1 \longrightarrow L \longrightarrow G \longrightarrow A \longrightarrow 1, \]
where $L$ is a connected linear algebraic group, and  $A$ an abelian variety.
Consider the unipotent radical
$R_u(L)$; this is a closed connected normal subgroup of $G$.
If $R_u(L)$ is non-trivial, consider the last non-trivial
term of its lower central series, $U$. Then $U$ is a closed
connected normal subgroup of $G$, isomorphic to the additive
group of a finite-dimensional vector space. Hence, $U$ is not
contained in $H$ by Claim \ref{c1}; but then $G$ has an additive one-parameter
subgroup with image not contained in $H$, contradicting Claim \ref{c2}.
So $R_u(L)$ is trivial, i.e., $L$ is reductive.
Thus, the derived subgroup $[L,L]$ is semi-simple,
and hence generated by images of additive
one-parameter subgroups. So $[L,L] \subseteq H$ by Claim \ref{c2}.
As $[L,L]$ is a normal subgroup of $G$, Claim \ref{c1} implies that $[L,L]$
is trivial. So $L$ is an algebraic torus. Hence $G$ is a semi-abelian
variety.

\par \vskip 1pc
This completes the proof of Theorem \ref{ThA}.

\par \vskip 1pc
{\it Next we prove Corollary \ref{CorA}.}

Denote by $N$ the kernel of the action of $G$ on $G/H$, i.e., the
largest (closed ) normal subgroup of $G$ contained in $H$.
Replacing $G,H$ with $G/N, H/N$, respectively, we may assume that
$G$ acts faithfully on $G/H$.

By \cite[Thm.~2]{Br10}, there exists a projective compactification
$X$ on $G/H$ such that the natural $G$-action on $G/H$ from the left
extends to an
algebraic action on $X$. Using equivariant resolution of singularities
(see \cite{Ko07}), we may assume that $X$ is smooth and
$X \setminus V$ is a simple normal crossing divisor. Then the desired
assertion follows from Theorem \ref{ThA}.
This proves Corollary \ref{CorA}.

\section{Concluding remarks}

\noindent
3.1. With the assumptions of Theorem \ref{ThA}, the action of the Lie
algebra of $G$ on $X$ preserves $D$, and hence yields a linear map
\[ \op_{X,D} : \Lie(G) \longrightarrow H^0(X,T_X(- \log D)), \]
where $T_X(- \log D)$ denotes the log tangent bundle. Moreover, the
pull-back of $T_X(- \log D)$ to $V$ is just the tangent bundle $T_V$,
which is generated by the image of $\op_{X,D}$. It follows that
the induced map
\[ \wedge^n \Lie(G) \longrightarrow H^0(X, -(K_X + D)) \]
is nonzero, where $n := \dim(V) = \dim(X)$. In particular,
$-(K_X + D)$ is effective. Thus, the pseudo-effectivity
of $K_X + D$ in Theorem \ref{ThA} (1) immediately implies
Theorem \ref{ThA} (2): $K_X + D \sim 0$.

\medskip

\noindent
3.2.
Theorem \ref{ThA} is applicable to an arbitrary (not necessarily smooth)
equivariant compactification $Y \supseteq V$ as long as $Y$ is normal,
the boundary $E:= Y \setminus V$ is a divisor, and $K_Y + E$ is $\Q$-Cartier.
Indeed, by equivariant resolution of singularities, there exists a
$G$-equivariant birational morphism $f : X \to Y$ such that $X$ is smooth,
projective and equipped with a $G$-action, and $D := X \setminus V$
(which is the inverse of $E$) is a simple normal crossing divisor.
By the logarithmic ramification divisor formula, we have
$$K_X + D = f^*(K_Y + E) + R_f$$
where $R_f$ is an effective $f$-exceptional
divisor. Thus the pairs $(X, D)$ and $(Y, E)$
(as compactifications of $V$) have the same log Kodaira dimensions:
$\kappa(X, K_X + D) = \kappa(Y, K_Y + E)$,
and the log canonical divisor $K_X + D$ is pseudo-effective if and only if
so is $K_Y + E$.

In particular, if $K_Y + E$ is pseudo-effective, then so is $K_X + D$,
and hence both $V \subseteq X$ and $V \subseteq Y$ are the embeddings
of the semi-abelian variety $V$ in $X$ and $Y$, respectively.

\medskip

\noindent
3.3. The normal projective equivariant compactifications of a
given semi-abelian variety are well-understood by work of
Alexeev (see \cite{Al02}), where they are called semi-abelic varieties.
In particular, for a semi-abelian variety $G$ given by an extension
$$ 1 \to T \to G \to A \to 1$$
where $T$ is an algebraic torus and $A$ an abelian
variety, and a normal projective equivariant compactification $Y$ of $G$,
the morphism $q : G \to A$ extends to a $T$-invariant morphism
\[ f : Y \longrightarrow A, \]
which is a fibration in projective toric varieties with the big torus $T$.
Thus, $f$ is the Albanese map of the variety $Y$.

\end{document}